\documentclass{amsart}
\usepackage{amsthm,amsfonts,amsmath,graphicx}
\usepackage[numbers,sort&compress]{natbib}
\usepackage[lmargin=4cm,rmargin=4cm,tmargin=3cm,bmargin=28mm]{geometry}
\renewcommand{\baselinestretch}{1.175}
\theoremstyle{plain}
\newtheorem{theorem}{Theorem}
\newtheorem{lemma}{Lemma}

\newcommand{\seclabel}[1]{\label{sec:#1}}

\newcommand{\HALF}{\ensuremath{\protect\dfrac{1}{2}}}

\newcommand{\n}[1]{\ensuremath{|#1|}}
\newcommand{\nG}{\n{G}}

\begin{document}

\title[Independent Sets in Graphs with an Excluded Clique Minor]{\textnormal{\footnotesize Note}\\ 
Independent Sets in Graphs\\ with an Excluded Clique Minor}

\author{David R. Wood}

\address{Departament de Matem{\`a}tica Aplicada II, Universitat Polit{\`e}cnica de Catalunya, Barcelona, Spain}

\thanks{Supported by the Government of Spain grant MEC SB2003-0270 and by the projects MCYT-FEDER BFM2003-00368 and Gen.\ Cat 2001SGR00224.}

\email{david.wood@upc.edu}

\date{\today}


\begin{abstract} 
Let $G$ be a graph with $n$ vertices, with independence number $\alpha$, and with with no $K_{t+1}$-minor for some $t\geq5$. It is proved that $(2\alpha-1)(2t-5)\geq2n-5$.
\end{abstract}

\maketitle

%%%%%%%%%%%%%%%%%%%%%%%%%%%%%%%%%%%%%%%%%%%%%%%%%%%%%%%%%%%%%%%%%%%%%%%%%%%
\section{Introduction}
\seclabel{Introduction}
%%%%%%%%%%%%%%%%%%%%%%%%%%%%%%%%%%%%%%%%%%%%%%%%%%%%%%%%%%%%%%%%%%%%%%%%%%%

In 1943, \citet{Hadwiger43} made the following conjecture, which is widely considered to be one of the most important open problems in graph theory\footnote{All graphs considered in this note are undirected, simple and finite. Let $G$ be a graph with vertex set $V(G)$. Let $\nG:=|V(G)|$. Let $X\subseteq V(G)$. $X$ is \emph{connected} if the subgraph of $G$ induced by $X$ is connected. $X$ is \emph{dominating} if every vertex of $G\setminus X$ has a neighbour in $X$. $X$ is \emph{independent} if no two vertices in $X$ are adjacent. The \emph{independence number} $\alpha(G)$ is the maximum cardinality of an independent set of $G$. $X$ is a \emph{clique} if every pair of vertices in $X$ are adjacent. The \emph{clique number} $\omega(G)$ is the maximum cardinality of a clique in $G$. A \emph{$k$-colouring} of $G$ is a function that assigns one of $k$ colours to each vertex of $G$ such that adjacent vertices receive distinct colours. The \emph{chromatic number} $\chi(G)$ is the minimum integer $k$ such that $G$ is $k$-colourable. A \emph{minor} of $G$ is a graph that can be obtained from a subgraph of $G$ by contracting edges. The \emph{Hadwiger number} $\eta(G)$ is the maximum integer $n$ such that the complete graph $K_n$ is a minor of $G$.}; see \citep{Toft-HadwigerSurvey96} for a survey.

\medskip\noindent\textbf{Hadwiger's Conjecture.} For every integer $t\geq1$, every graph with no $K_{t+1}$-minor is $t$-colourable. That is, $\chi(G)\leq\eta(G)$ for every graph $G$.
\medskip

Hadwiger's Conjecture is trivial for $t\leq2$, and is straightforward for 
$t=3$; see \citep{Hadwiger43,Dirac52,Woodall-JGT92}. In the cases $t=4$ and $t=5$, \citet{Wagner37}  and \citet{RST-Comb93} respectively proved that Hadwiger's Conjecture is equivalent to the Four-Colour Theorem \citep{AH77_I,AH77_II,RSST97}. Hadwiger's Conjecture is open for all $t\geq6$. 
Progress on the $t=6$ case has been recently been obtained by \citet{KT-Comb05} (without using the Four-Colour Theorem). The best known upper bound is $\chi(G)\leq c\cdot\eta(G)\sqrt{\log\eta(G)}$ for some constant $c$, independently due to \citet{Kostochka82} and \citet{Thomason84,Thomason01}. 

%In particular, it is unknown whether $\chi(G)\leq c\cdot\eta(G)$ for some constant $c$. In this note we propose a conjecture that would imply that $\chi(G)\leq 2\cdot\eta(G)$. 

Since $\alpha(G)\cdot\chi(G)\geq\nG$ for every graph $G$, Hadwiger's Conjecture implies that 
\begin{equation}
\label{WeakHadwiger}
\alpha(G)\cdot\eta(G)\geq\nG,
\end{equation}
as observed by \citet{Woodall-JGT87}. In general, \eqref{WeakHadwiger} is weaker than Hadwiger's Conjecture, but if $\alpha(G)=2$, then \citet{PST-DMGT03} proved that \eqref{WeakHadwiger} is in fact equivalent to Hadwiger's Conjecture. The first significant progress towards \eqref{WeakHadwiger} was made by \citet{DM82} (also see \citep{MM-DM87}), who proved that 
\begin{equation}
\label{DM}
(2\alpha(G)-1)\cdot\eta(G)\geq\nG\enspace.
\end{equation}
This result was improved by \citet{KPT-JCTB05} to
\begin{equation}
\label{KPT1}
(2\alpha(G)-1)\cdot\eta(G)\geq\nG+\omega(G)\enspace.
\end{equation}
Assuming $\alpha(G)\geq3$, \citet{KPT-JCTB05} proved that 
\begin{equation}
\label{KPT}
(4\alpha(G)-3)\cdot\eta(G)\geq2\nG,
\end{equation}
which was further improved by \citet{KS05} to 
\begin{equation}
\label{KS}
(2\alpha(G)-2)\cdot\eta(G)\geq\nG.
\end{equation}

The following theorem is the main contribution of this note.

\begin{theorem}
\label{MainTheorem}
Every graph $G$ with $\eta(G)\geq5$ satisfies
\begin{equation*}
(2\alpha(G)-1)(2\eta(G)-5)\geq2\nG-5\enspace.
\end{equation*}
\end{theorem}

Observe that Theorem~\ref{MainTheorem} represents an improvement over \eqref{DM}, \eqref{KPT} and \eqref{KS} whenever $\eta(G)\geq5$ and $\nG\geq\tfrac{2}{5}\eta(G)^2$. For example, Theorem~\ref{MainTheorem} implies that $\alpha(G)>\frac{\nG}{7}$ for every graph $G$ with $\eta(G)\leq 6$, whereas each of \eqref{DM}, \eqref{KPT} and \eqref{KS} imply that $\alpha(G)>\frac{\nG}{12}$.

\section{Proof of Theorem~\ref{MainTheorem}}
%%%%%%%%%%%%%%%%%%%%%%%%%%%%%%%%%%%%%%%%%%%%%%%%%%%%%%%%%%%%%%%%%%%%%%%%%%%%

Theorem~\ref{MainTheorem} employs the following lemma by \citet{DM82}. The proof is included for completeness. 

\begin{lemma}[\citep{DM82}]
\label{DominatingSet}
Every connected graph $G$ has a connected dominating set $D$ and an independent set $S\subseteq D$ such that $|D|=2|S|-1$.
\end{lemma}

\begin{proof}
Let $D$ be a maximal connected set of vertices of $G$ such that $D$ contains an independent set $S$ of $G$ and $|D|=2|S|-1$. There is such a set since $D:=S:=\{v\}$ satisfies these conditions for each vertex $v$. We claim that $D$ is dominating. Otherwise, since $G$ is connected, there is a vertex $v$ at distance two from $D$, and there is a neighbour $w$ of $v$ at distance one from $D$. Let  $D:=D\cup\{v,w\}$ and $S':=S\cup\{v\}$. Thus $D'$ is connected and contains an independent set $S'$ such that $|D'|=2|S'|-1$. Hence $D$ is not maximal. This contradiction proves that $D$ is dominating. 
\end{proof}

The next lemma is the key to the proof of Theorem~\ref{MainTheorem}.

\begin{lemma}
\label{MainLemma}
Suppose that for some integer $t\geq1$ and for some real number $p\geq t$, every graph $G$ with $\eta(G)\leq t$ satisfies 
$p\cdot\alpha(G)\geq\nG$. Then every graph $G$ with $\eta(G)\geq t$ satisfies 
\begin{equation*}
\alpha(G)\geq\frac{2\nG-p}{4\eta(G)+2p-4t}+\HALF\enspace.
\end{equation*}
\end{lemma}

\begin{proof}
We proceed by induction on $\eta(G)-t$. If $\eta(G)=t$ the result holds by assumption. Let $G$ be a graph with $\eta(G)>t$. We can assume that $G$ is connected. By Lemma~\ref{DominatingSet}, $G$ has a connected dominating set $D$ and an independent set $S\subseteq D$ such that $|D|=2|S|-1$. Now $\alpha(G)\geq|S|=\frac{|D|+1}{2}$. Thus we are done if 
\begin{equation}
\label{Dbig}
\frac{|D|+1}{2}\geq\frac{2\nG-p}{4\eta(G)+2p-4t}+\HALF\enspace.
\end{equation}
Now assume that \eqref{Dbig} does not hold. That is,
\begin{equation*}
|D|\leq\frac{2\nG-p}{2\eta(G)+p-2t}\enspace.
\end{equation*}
Thus
\begin{equation*}
\n{G\setminus D}=\nG-|D|
\geq\frac{(2\eta(G)+p-2t-2)\nG+p}{2\eta(G)+p-2t}\enspace.
\end{equation*}
Since $D$ is dominating and connected, $\eta(G\setminus D)\leq\eta(G)-1$. 
Thus by induction,
\begin{align*}
\alpha(G)\geq\alpha(G\setminus D)
&\geq\frac{2\n{G\setminus D}-p}{4\eta(G\setminus D)+2p-4t}+\HALF\\
&\geq\frac{2(2\eta(G)+p-2t-2)\nG+2p}{(2\eta(G)+p-2t)(4\eta(G)-4+2p-4t)}-\frac{p}{4\eta(G)-4+2p-4t}+\HALF\\
%&=\frac{\nG}{2\eta(G)+p-2t}+\frac{2p}{(2\eta(G)+p-2t)(4\eta(G)-4+2p-4t)}-\frac{p}{4\eta(G)-4+2p-4t}+\HALF\\
%&=\frac{\nG}{2\eta(G)+p-2t}-\frac{p}{4\eta(G)+2p-4t}+\HALF\\
&=\frac{2\nG-p}{4\eta(G)+2p-4t}+\HALF\enspace.
\end{align*}
This completes the proof.
\end{proof}

\begin{lemma}
\label{SecondLemma}
Suppose that Hadwiger's Conjecture is true for some integer $t$. 
Then every graph $G$ with $\eta(G)\geq t$ satisfies 
\begin{equation*}
(2\eta(G)-t)(2\alpha(G)-1)\geq 2\nG-t\enspace.
\end{equation*}
\end{lemma}

\begin{proof}
If Hadwiger's Conjecture is true for $t$ then $t\cdot\alpha(G)\geq\nG$ for every graph $G$ with $\eta(G)\leq t$. Thus Lemma~\ref{MainLemma} with $p=t$ implies that every graph $G$ with $\eta(G)\geq t$ satisfies 
\begin{equation*}
\alpha(G)\geq\frac{2\nG-t}{4\eta(G)-2t}+\HALF\enspace,
\end{equation*}
which implies the result. 
\end{proof}

Theorem~\ref{MainTheorem} follows from Lemma~\ref{SecondLemma} with $t=5$ since Hadwiger's Conjecture holds for $t=5$ \citep{RST-Comb93}.

\section{Concluding Remarks}
%%%%%%%%%%%%%%%%%%%%%%%%%%%%%%%%%%%%%%%%%%%%%%%%%%%%%%%%%%%%%%%%%%%%%%%%%%%%

The proof of Theorem~\ref{MainTheorem} is substantially simpler than the proofs of \eqref{KPT1}--\eqref{KS}, ignoring its dependence on the proof of Hadwiger's Conjecture with $t=5$, which in turn is based on the four-colour theorem. A bound that still improves upon \eqref{DM}, \eqref{KPT} and \eqref{KS} but with a completely straightforward proof is obtained from Lemma~\ref{SecondLemma} with $s=3$: Every graph $G$ with $\eta(G)\geq3$ satisfies 
$(2\eta(G)-3)(2\alpha(G)-1)\geq2\nG-3$.

We finish with an open problem. The method of \citet{DM82} was generalised by \citet{ReedSeymour-JCTB98} to prove that the fractional chromatic number $\chi_f(G)\leq 2\eta(G)$. For sufficiently large $\eta(G)$, is $\chi_f(G)\leq 2\eta(G)-c$ for some constant $c\geq1$?

%%%%%%%%%%%%%%%%%%%%%%%%%%%%%%%%%%%%%%%%%%%%%%%%%%%%%%%%%%%%%%%%%%%%%%%%%%%%%%%
%\bibliographystyle{myNatbibStyle}
%\bibliography{myBibliography,myConferences}

\begin{thebibliography}{20}
\providecommand{\natexlab}[1]{#1}
\providecommand{\url}[1]{\texttt{#1}}
\providecommand{\urlprefix}{}
\expandafter\ifx\csname urlstyle\endcsname\relax
  \providecommand{\doi}[1]{doi:\discretionary{}{}{}#1}\else
  \providecommand{\doi}{doi:\discretionary{}{}{}\begingroup
  \urlstyle{rm}\Url}\fi

\bibitem[{Appel and Haken(1977)}]{AH77_I}
\textsc{Kenneth Appel and Wolfgang Haken}.
\newblock Every planar map is four colorable. {I}. {D}ischarging.
\newblock \emph{Illinois J. Math.}, 21(3):429--490, 1977.

\bibitem[{Appel et~al.(1977)Appel, Haken, and Koch}]{AH77_II}
\textsc{Kenneth Appel, Wolfgang Haken, and John Koch}.
\newblock Every planar map is four colorable. {I}{I}. {R}educibility.
\newblock \emph{Illinois J. Math.}, 21(3):491--567, 1977.

\bibitem[{Dirac(1952)}]{Dirac52}
\textsc{Gabriel~A. Dirac}.
\newblock A property of {$4$}-chromatic graphs and some remarks on critical
  graphs.
\newblock \emph{J. London Math. Soc.}, 27:85--92, 1952.

\bibitem[{Duchet and Meyniel(1982)}]{DM82}
\textsc{Pierre Duchet and Henri Meyniel}.
\newblock On {H}adwiger's number and the stability number.
\newblock \emph{Annals of Discrete Mathematics}, 13:71--73,  1982.

\bibitem[{Hadwiger(1943)}]{Hadwiger43}
\textsc{Hugo Hadwiger}.
\newblock \"{U}ber eine {K}lassifikation der {S}treckenkomplexe.
\newblock \emph{Vierteljschr. Naturforsch. Ges. Z\"urich}, 88:133--142, 1943.

\bibitem[{Kawarabayashi et~al.(2005)Kawarabayashi, Plummer, and
  Toft}]{KPT-JCTB05}
\textsc{Ken-ichi Kawarabayashi, Michael~D. Plummer, and Bjarne Toft}.
\newblock Improvements of the theorem of {D}uchet and {M}eyniel on {H}adwiger's
  conjecture.
\newblock \emph{J. Combin. Theory Ser. B}, 95(1):152--167, 2005.

\bibitem[{Kawarabayashi and Song(2005)}]{KS05}
\textsc{Ken-ichi Kawarabayashi and Zi-Xia Song}.
\newblock Independence number and clique minors, 2005.
\newblock \urlprefix\url{http://www.math.ohio-state.edu/\~{}song/kensong.pdf}.

\bibitem[{Kawarabayashi and Toft(2005)}]{KT-Comb05}
\textsc{Ken-ichi Kawarabayashi and Bjarne Toft}.
\newblock Any 7-chromatic graph has {$K\sb 7$} or {$K\sb {4,4}$} as a minor.
\newblock \emph{Combinatorica}, 25(3):327--353, 2005.

\bibitem[{Kostochka(1982)}]{Kostochka82}
\textsc{Alexandr~V. Kostochka}.
\newblock The minimum {H}adwiger number for graphs with a given mean degree of
  vertices.
\newblock \emph{Metody Diskret. Analiz.}, 38:37--58, 1982.

\bibitem[{Maffray and Meyniel(1987)}]{MM-DM87}
\textsc{Fr{\'e}d{\'e}ric Maffray and Henri Meyniel}.
\newblock On a relationship between {H}adwiger and stability numbers.
\newblock \emph{Discrete Math.}, 64(1):39--42, 1987.

\bibitem[{Plummer et~al.(2003)Plummer, Stiebitz, and Toft}]{PST-DMGT03}
\textsc{Michael~D. Plummer, Michael Stiebitz, and Bjarne Toft}.
\newblock On a special case of {H}adwiger's conjecture.
\newblock \emph{Discuss. Math. Graph Theory}, 23(2):333--363, 2003.

\bibitem[{Reed and Seymour(1998)}]{ReedSeymour-JCTB98}
\textsc{Bruce Reed and Paul Seymour}.
\newblock Fractional colouring and {H}adwiger's conjecture.
\newblock \emph{J. Combin. Theory Ser. B}, 74(2):147--152, 1998.

\bibitem[{Robertson et~al.(1997)Robertson, Sanders, Seymour, and
  Thomas}]{RSST97}
\textsc{Neil Robertson, Daniel~P. Sanders, Paul Seymour, and Robin Thomas}.
\newblock The four-colour theorem.
\newblock \emph{J. Combin. Theory Ser. B}, 70(1):2--44, 1997.

\bibitem[{Robertson et~al.(1993)Robertson, Seymour, and Thomas}]{RST-Comb93}
\textsc{Neil Robertson, Paul Seymour, and Robin Thomas}.
\newblock Hadwiger's conjecture for ${K}\sb 6$-free graphs.
\newblock \emph{Combinatorica}, 13(3):279--361, 1993.

\bibitem[{Thomason(1984)}]{Thomason84}
\textsc{Andrew Thomason}.
\newblock An extremal function for contractions of graphs.
\newblock \emph{Math. Proc. Cambridge Philos. Soc.}, 95(2):261--265, 1984.

\bibitem[{Thomason(2001)}]{Thomason01}
\textsc{Andrew Thomason}.
\newblock The extremal function for complete minors.
\newblock \emph{J. Combin. Theory Ser. B}, 81(2):318--338, 2001.

\bibitem[{Toft(1996)}]{Toft-HadwigerSurvey96}
\textsc{Bjarne Toft}.
\newblock A survey of {H}adwiger's conjecture.
\newblock \emph{Congr. Numer.}, 115:249--283, 1996.

\bibitem[{Wagner(1937)}]{Wagner37}
\textsc{Klaus Wagner}.
\newblock {\"U}ber eine {E}igenschaft der ebene {K}omplexe.
\newblock \emph{Math. Ann.}, 114:570--590, 1937.

\bibitem[{Woodall(1987)}]{Woodall-JGT87}
\textsc{Douglas~R. Woodall}.
\newblock Subcontraction-equivalence and {H}adwiger's conjecture.
\newblock \emph{J. Graph Theory}, 11(2):197--204, 1987.

\bibitem[{Woodall(1992)}]{Woodall-JGT92}
\textsc{Douglas~R. Woodall}.
\newblock A short proof of a theorem of {D}irac's about {H}adwiger's
  conjecture.
\newblock \emph{J. Graph Theory}, 16(1):79--80, 1992.

\end{thebibliography}
%%%%%%%%%%%%%%%%%%%%%%%%%%%%%%%%%%%%%%%%%%%%%%%%%%%%%%%%%%%%%%%%%%%%%%%%%%%%%%%

\def\soft#1{\leavevmode\setbox0=\hbox{h}\dimen7=\ht0\advance \dimen7
  by-1ex\relax\if t#1\relax\rlap{\raise.6\dimen7
  \hbox{\kern.3ex\char'47}}#1\relax\else\if T#1\relax
  \rlap{\raise.5\dimen7\hbox{\kern1.3ex\char'47}}#1\relax \else\if
  d#1\relax\rlap{\raise.5\dimen7\hbox{\kern.9ex \char'47}}#1\relax\else\if
  D#1\relax\rlap{\raise.5\dimen7 \hbox{\kern1.4ex\char'47}}#1\relax\else\if
  l#1\relax \rlap{\raise.5\dimen7\hbox{\kern.4ex\char'47}}#1\relax \else\if
  L#1\relax\rlap{\raise.5\dimen7\hbox{\kern.7ex
  \char'47}}#1\relax\else\message{accent \string\soft \space #1 not
  defined!}#1\relax\fi\fi\fi\fi\fi\fi} \def\cprime{$'$}

\end{document}